\documentclass{amsart}
\usepackage{amssymb}
\usepackage{latexsym}
\begin{document}
\setlength{\parskip}{.15 in}
\setlength{\parindent}{.25 in}
\newtheorem{lemma}{Lemma}
\newtheorem{theorem}{Theorem}
\newtheorem{corollary}{Corollary}
\newtheorem{example}{Example}
\newtheorem{remark}{Remark}
\newtheorem*{guess}{Conjecture}

\title{Cluster points and asymptotic values
of planar harmonic functions}
\author{Genevra Neumann}
\address{Department of Mathematics,
Kansas State University,
Manhattan, KS 66506-2602}
\email{neumann@math.ksu.edu}
\date{August 30, 2005}
\subjclass[2000]{Primary 30C99, 26B99;
Secondary 31A05, 26C99}
\keywords{Planar harmonic functions,
$C^1$ functions in $\mathbb{R}^2$,
cluster set, asymptotic values}

\begin{abstract}
A sufficient condition for a cluster point 
of a planar harmonic function to be
an asymptotic value is given,
based on a partitioning into regions of constant 
valence.
A sufficient condition for the cluster set of a planar
harmonic function to have non-empty interior is
given.
An example is given of a planar harmonic function
where the image of the critical set is not closed
and such that the cluster set has non-empty interior
and is a proper subset of the image.
\end{abstract}

\maketitle

\section{Introduction}
The behavior of complex-valued harmonic functions
in the plane was studied in \cite{neu:valence}.
Many of these results also hold for
$C^1$ functions $f(z) = (u(z), v(z))$ in an open set
$R \subseteq \mathbb{R}^2$,
identifying $z = (x,y) \in \mathbb{R}^2$
with $z = x + iy \in \mathbb{C}$ 
and $f$ with $u + iv$.
Recall that
critical set of $f$ is the set of points where
the Jacobian of $f$ vanishes.
The cluster set of $f$ consists of those finite values
$w$ such that there exists a sequence 
$\{z_n\} \subset R$ such that 
$z_n \rightarrow \partial R \cup \{\infty\}$
and $f(z_n) \rightarrow w$.
It was shown that the image of the critical set
combined with the cluster set partitions
$\mathbb{R}^2$ into regions where each point has
the same number of distinct preimages.
This result and some of its consequences are
summarized in Theorem \ref{thm:part} below.
These results are nicest when the partitioning set
has empty interior
and that occurs when the cluster set has empty
interior.

It was shown that if the cluster set has
non-empty interior,
then points with infinite valence are dense in the
interior of the cluster set.
Harmonic functions will be assumed to be complex-valued
unless stated otherwise.
If $f$ is also assumed to be harmonic in $\mathbb{C}$,
it was shown that if the cluster set of 
$f = h + \overline{g}$ 
(where $h$ and $g$ are entire)
has empty interior,
then both $g$ and $h$ are either polynomials 
or entire transcendental functions.
Theorem \ref{thm:clusvals} below
gives a sufficient condition for the critical
set to have non-empty interior.

Many examples were given in \cite{neu:valence}.
In these examples,
the image of the critical set was always closed.
Further,
when a cluster set had non-empty interior,
it filled the entire plane.
Is the image of the critical set always closed?
Can a cluster set with non-empty interior be a proper
subset of the image?
Example \ref{ex:nono} below shows that the image of the
critical set is not always closed and that a cluster
set with non-empty interior can be a proper subset of
the image.

In the examples of entire harmonic functions in 
Section 4 of \cite{neu:valence},
where the cluster set has empty interior,
the sequence given for $\{z_n\}$
where $z_n \rightarrow \infty$ and where
$f(z_n)$ approaches a finite limit,
can be extended to a curve approaching $\infty$
whose image approaches the cluster value.
In these examples,
the cluster set has empty interior and
each finite cluster point is an asymptotic value.
Does this hold in general for
entire harmonic functions?
Note that entire harmonic functions can exclude
open regions of the plane.

Consider the following version of Iversen's theorem:

\begin{theorem}[\cite{car:integral} p.25]
Suppose that $f(z)$ is meromorphic for $|z| > \ell$
and has an essential singularity at infinity.
Then if $f(z)$ has a Picard exceptional value $a$
(\textit{e.g.},
$f(z) - a$ has a finite number of zeros),
it is an asymptotic value of $f(z)$.
In particular,
infinity is an asymptotic value of $f$.
\end{theorem}

\noindent
Recall that a sense-preserving harmonic function 
$f(z)$ is said to be
K-quasiregular in a domain
$U \subseteq \mathbb{R}^2$ 
if there exists a constant $K \geq 1$
such that
$(|f_z| + |f_{\overline{z}}|)^2 \leq K\, J_f$
holds for all $z \in U$
(see \cite{cgh:bloch}).
Note that the term K-quasiregular is defined for
functions 
defined in subsets of $\mathbb{R}^n$ for some
$n \geq 2$ and does not require that
$f$ be harmonic;
see \cite{ric:book}.
K-quasiregular functions can be thought of as
analogues to analytic functions and
K-quasimeromorphic functions as 
analogues to meromorphic functions.
O. Martio, S. Rickman, and J. V\"{a}is\"{a}l\"{a} 
showed the following analogue of Iversen's theorem
for quasimeromorphic functions:

\begin{theorem}[\cite{ric:book}, page 170]
Let $f: G \rightarrow \overline{\mathbb{R}}^n$
be a K-quasimeromorphic mapping and let 
$b \in \partial G$ be an isolated essential singularity
of $f$.
Then every point in 
$\overline{\mathbb{R}}^n \backslash fG$
is an asymptotic value of $f$.
\end{theorem}

\noindent
D. A. Brannan 
\cite{br:infty}
considered the case
$f:\mathbb{C} \rightarrow \mathbb{C}$
where $f$ is continuous:

\begin{theorem}[Brannan \cite{br:infty}]
Let $f:\mathbb{C} \rightarrow \mathbb{C}$
be continuous.
Then
at least one of the following occurs:
(i) $f$ has $\infty$ as an asymptotic value,
(ii) $f$ is bounded on a path going to $\infty$,
or (iii)
$f$ is uniformly bounded on a sequence of closed curves
which surround the origin and whose distance
from the origin approaches $\infty$.
\end{theorem}

Is there an analog to Iversen's theorem for
planar harmonic functions?
It seems reasonable to expect that
the behavior of $f$ harmonic in $\mathbb{C}$
to be nicer than when $f$ is only continuous,
but perhaps not as nice as in the K-quasimeromorphic
and meromorphic cases.
Clearly,
if $\lim_{z \rightarrow \infty} f(z) = \infty$,
then $\infty$ is an asymptotic value of $f$.
Based on the examples in \cite{neu:valence},
one might expect

\begin{guess}
Let $f:\mathbb{C} \rightarrow \mathbb{C}$
be harmonic.
If the critical set of $f$ is nowhere dense
and if the cluster set of $f$ has empty interior,
then each cluster point of $f$ is an
asymptotic value.
In particular,
if each value in $f(\mathbb{C})$ is taken at most
a finite number of times,
then every cluster point of $f$ is an asymptotic
value.
\end{guess}

\noindent
Theorem \ref{thm:locconpath} gives a sufficient
condition for a cluster value of an
entire $C^1$ function on $\mathbb{R}^2$ 
to be an asymptotic value.
Given that this is a $C^1$ result,
it is not surprising that it
does not settle this conjecture.
The proof uses a path-lifting argument.
It starts with a sequence $\{z_n\}$ such that
$z_n \rightarrow \infty$ and 
$w_n = f(z_n)$ approaches the cluster point $w_0$.
It builds a simple asymptotic polygonal path that
contains a subsequence of the $\{w_n\}$
and approaches $w_0$.
Theorem \ref{thm:part} is then used to lift this
path and it is shown that one of the lifted paths
approaches $\infty$ asymptotically.

\section{Notation and Main Results}

\noindent
Let $R \subseteq \mathbb{R}^2$ be open and let 
$f: R \rightarrow \mathbb{R}^2$ be $C^1$.
Let

\noindent
$B(z_0, r) = \{z: |z - z_0| < r\}$\\
$J_f = |f_z|^2 - |f_{\overline{z}}|^2$\\
$S = \{z \in R : J_f(z) = 0\}$\\
$C(f) = \{w \in \mathbb{R}^2: \exists\, 
\{z_n\} \subset R$ where
$z_n \rightarrow z_0 \in 
\partial R \cup \{\infty\}$ and
$f(z_n) \rightarrow w\}$\\
$Val(f, w) = \#\{z \in R : f(z) = w\}$\\
$\Lambda_f = |f_z| + |f_{\overline{z}}|$\\
$\lambda_f = ||f_z| - |f_{\overline{z}}||$\\
$\mu_f = f_{\overline{z}} / f_z$

\noindent
$S$ is called the critical set of $f$ and $C(f)$ the 
cluster set.
It is well known that $C(f)$ is closed.
We also note that $f(S) \cup C(f)$ is closed
(see \cite{neu:valence}).
If $R = \mathbb{R}^2$, we denote the cluster set by 
$C(f,\infty)$.
$Val(f, w)$ is the valence of $f$ at $w$ and counts
the number of distinct preimages of $w$ in $R$;
it does not count multiplicity.
$\mu_f$ is referred to as the complex dilatation
of $f$.

We summarize some results from
\cite{neu:valence} as follows:

\begin{theorem}
\label{thm:part}
Let $R \subseteq \mathbb{R}^2$ be open and
let $f:R \rightarrow \mathbb{R}^2$ be $C^1$.
Then 
$f(S) \cup C(f)$ partitions $\mathbb{R}^2$
into regions of constant, finite valence.
Moreover,
each component 
$R_0 \subseteq R \backslash f^{-1}(f(S) \cup C(f))$
is a covering space for $f(R_0)$ with 
$f|_{R_0}$ as the covering map.
Note that $f(R_0)$ is a component of 
$\mathbb{R}^2 \backslash (f(S) \cup C(f))$.
\end{theorem}

We will refer to $f(S) \cup C(f)$ as the partitioning
set of the image and
to $f^{-1}(f(S) \cup C(f))$ as the partitioning set of 
the preimage.
Unless explicitly stated otherwise,
a component of the image is a component of
$\mathbb{R}^2 \backslash (f(S) \cup C(f))$ and
a component of the preimage is a component of
$R \backslash f^{-1}(f(S) \cup C(f))$.

Let $w_0 \in C(f)$.
If there exists an asymptotic path
$\gamma \rightarrow z_0 \in \partial R \cup \{\infty\}$
such that $f(\gamma) \rightarrow w_0$,
then $w_0$ is said to be an asymptotic value 
of $f$.
Recall that an asymptotic path is
a curve which is the continuous
image of $[0, \infty)$ 
such that $\gamma(t) \in R$ for all $0 \leq t < \infty$
and 
$\gamma \rightarrow z_0  
\in \partial R \cup \{\infty\}$.
If $\gamma \rightarrow z_0 \in \partial R$ 
where $z_0$ is finite,
we also refer to $\gamma$ as an end-cut.

Recall that a set $U$ is uniformly locally connected
if for each $\epsilon > 0$,
there exists $\delta > 0$ such that if
$z_1, z_2 \in U$ with $|z_1 - z_2| < \delta$,
then there exists a connected subset of $U$ 
joining $z_1$
and $z_2$ of diameter less than $\epsilon$
(\cite{new:plane}, page 160).
If instead of requiring a connected subset joining the
two points,
we require an arc,
the set is set to be uniformly locally arcwise 
connected
(ulac)
(\cite{hy:top}, page 129).

Our two main results are: 

\begin{theorem}
\label{thm:clusvals}
Let $f:\mathbb{C} \rightarrow \mathbb{C}$ be harmonic.
Suppose that $\{z_n\}$ is such that 
$z_n \rightarrow \infty$,
$w_n = f(z_n) \rightarrow w_0$,
and $J_f(z_n) > 0$ for all $n$.
Suppose that the following hold:
\begin{enumerate}
\item
$d(\{z_n\}, S) = \delta > 0$,
\item
$\inf_n J_f(z_n) = \eta^2 > 0$,
and
\item
There exists $\rho$ in $(0, \delta)$
such that
$\sup_{z \in B(z_n, \rho)_{n=1}^\infty} 
J_f(z) / J_f(z_n) < \infty$.
\end{enumerate}
If $\sup_n \{|\mu_f(z)|: z \in 
\overline{B(z_n, \rho)}\} < 1$,
then $int\ C(f, \infty) \neq \varnothing$
and $Val(f, w) = \infty$ in a neighborhood of $w_0$.
\end{theorem}

\begin{remark}
\textnormal{
The condition that $J_f(z_n) > 0$ is not needed in
Theorem \ref{thm:clusvals},
given that the sequence stays off the critical set.
Given $\{z_n\}$ such that $z_n \rightarrow \infty$ and
$d(\{z_n\}, S) > 0$,
we can find a subsequence such that
$J_f(z_n)$ does not change signs;
this follows from noting that our sequence is off the critical set,
applying the pigeonhole principle, and passing to a
subsequence.
If $J_f(z_n) < 0$ for a subsequence,
we can work with $\overline{f}$ instead.}
\end{remark}

\begin{theorem}
\label{thm:locconpath}
Let $f$ be a $C^1$ function on $\mathbb{R}^2$ and 
suppose that 
$S$ is nowhere dense.
Let $w_0 \in C(f, \infty)$ be such that
either (i) $Val(f, w_0) < \infty$ or 
(ii) $Val(f, w_0) = \infty$ and $w_0 \notin f(S)$. 
Suppose that there exists $\epsilon > 0$ such that:
\begin{enumerate}
\item
$C(f, \infty)$ is 
nowhere dense in $B(w_0, \epsilon)$.
\item
$B(w_0, \epsilon)
\backslash (f(S) \cup C(f, \infty))$ consists of 
finitely many components,
each of which is uniformly locally arcwise 
connected.
\end{enumerate}
Then $w_0$ is an asymptotic value of $f$.
\end{theorem}

\section{Proof of Theorem \ref{thm:clusvals}}

S. Bochner
\cite{bo:bloch}
proved a version of Bloch's theorem for
sense-preserving,
K-quasiregular harmonic mappings in $\mathbb{R}^n$
for $n \geq 2$.
H. Chen, P. M. Gauthier, and W. Hengartner
\cite{cgh:bloch}
studied the value of the Bloch constant for $n=2$.
For example,

\begin{theorem}[Chen, Gauthier, and Hengartner
\cite{cgh:bloch}]
Let $f$ be a K-quasiregular harmonic mapping of 
the unit disk $\mathbb{D}$ such that
$J_f(0) = 1$.
Then $f(\mathbb{D})$ contains a schlicht disk
of radius at least 
$r_1 = \frac{\pi}{8 \sqrt{2} \sqrt{K}\, (1 + 2K)}$.
\end{theorem}

\noindent
However,
they also note that in harmonic analogues of
Bloch's theorem,
the schlicht disk contained in the image of the 
unit disk is not necessarily centered at $f(0)$.
Because
we are interested in finding disks in the image
that overlap in a neighborhood of our cluster point,
we need a somewhat different result.
Instead of adapting Bochner's proof,
we will use another result in \cite{cgh:bloch}:

\begin{theorem}[Chen, Gauthier, and Hengartner
\cite{cgh:bloch}]
\label{thm:cgh}
Let $f$ be a harmonic mapping of the unit disk
$\mathbb{D}$
such that $f(0) = 0$,
$\lambda_f(0) = 1$,
and
$\Lambda_f(z) \leq \Lambda$ for all 
$z \in \mathbb{D}$.
Then,
$f$ is univalent on $B(0, \rho_0)$
with $\rho_0 = \frac{\pi}{4(1 + \Lambda)}$
and $f(B(0, \rho_0))$ contains a schlicht disk
$B(0, r_0)$ with 
$r_0 = \frac{1}{2} \rho_0 = 
\frac{\pi}{8 (1 + \Lambda)}$.
\end{theorem}

\noindent
\textbf{Proof of Theorem \ref{thm:clusvals}.}
Let $M = \sup_n\{|\mu_f(z)|: z \in 
\overline{B(z_n, \rho)}\}$.
By assumption,
$M < 1$.
Let $K = (1 + M) / (1 - M)$.
Note that on each $B(z_n, \rho)$,
we have
$\Lambda_f^2 / J_f = (1 + |\mu_f|) / (1 - |\mu_f|)
\leq (1 + M) / (1 - M) = K$,
so $f$ is $K$-quasiregular on each $B(z_n, \rho)$.
Thus,
on each $B(z_n, \rho)$,
$\lambda_f^2 \geq J_f / K$.
In particular, 
for each $n$,
$(\lambda_f(z_n))^2 \geq J_f(z_n) / K$.

Following Bochner's proof,
we define the complex-valued harmonic function
\begin{center}
$F_n(z) = (f(z_n + \rho z) - f(z_n)) /
(\rho\, \lambda_f(z_n))$
\end{center}
for $z \in \mathbb{D}$.
We now check that each $F_n$ satisfies
the conditions of Theorem \ref{thm:cgh}
on $\mathbb{D}$.
Note that $F_n(0) = 0$
and $J_{F_n}(z) = 
J_f(z_n + \rho z) / (\lambda_f(z_n))^2$.
Also,
$(\Lambda_{F_n}(z))^2 = 
(\Lambda_f(z_n +\rho z))^2/ (\lambda_f(z_n))^2
\leq K J_f(z_n + \rho z) / 
(\lambda_f(z_n))^2
\leq K^2 J_f(z_n + \rho z) / J_f(z_n)
\leq \Lambda$ for each $n$ by condition (3),
where $\Lambda$ is a positive constant that is
independent of $n$.
A calculation shows that
$\lambda_{F_n}(0) = 1$ for each $n$.

By Theorem \ref{thm:cgh},
$F_n(\mathbb{D}) \supseteq B(0, r_0)$ for each $n$
where $r_0 = \pi / (8 (1 + \Lambda))$.
Thus,
$f(B(z_n, \rho)) \supseteq 
B(w_n, r_0 \rho \lambda_f(z_n))
\supseteq
B(w_n, r_1)$
for each $n$,
where,
by condition (2), 
$r_1 = r_0 \rho \eta / \sqrt{K}$.
Choose $0 < \epsilon_0 < r_1 / 2$.
Choose $N > 0$ such that $|w_n - w_0| < \epsilon_0$ for all $n \geq N$.
Let $w \in B(w_0, \epsilon_0)$.
Then $|w - w_n| \leq |w - w_0| + |w_0 - w_n|
< 2 \epsilon_0$.
Thus, $w \in B(w_n, r_1)$ and 
there exists $\zeta_n \in B(z_n, \rho)$ such that 
$f(\zeta_n) = w$.
By construction, $\zeta_n \rightarrow \infty$.
Thus, $w \in C(f, \infty)$
and $int\ C(f, \infty) \neq \varnothing$.
Moreover,
$Val(f, w) = \infty$.
\qed

\section{Proof of Theorem \ref{thm:locconpath}}

Suppose that we have 
a sequence $\{z_n\}$ such that
$z_n \rightarrow \infty$ and
$w_n = f(z_n) \rightarrow w_0$.
We set things up so that we may choose $w_n$ off
of the partitioning set of the image.
We show that $w_0$ is an asymptotic value
of $f$ by constructing a path $\gamma$ containing
a subsequence of the $\{w_n\}$ so that 
$\gamma \rightarrow w_0$ and then lifting
$\gamma$ to a component of the
preimage.
We then show that a lift of $\gamma$
asymptotically approaches $\infty$.
To simplify some of the arguments involving
the lifts of $\gamma$,
we will require that $\gamma$ be a simple path;
this construction is carried out by inductively
constructing a sequence of simple polygonal paths.

We begin with some elementary results
concering polygonal paths.
The proofs of Lemmas \ref{lem:simplepoly}
and \ref{lem:arctopoly} are sufficiently
simple that they will be omitted.

\begin{lemma}
\label{lem:simplepoly}
Let $\Gamma$ be a polygonal path consisting of
a finite number of line segments where
$\Gamma(0) = a$ and $\Gamma(1) = b$.
Then there exists a simple polygonal path
(consisting of a finite number of line segments) 
$\gamma \subseteq \Gamma$ such that
$\gamma(0) = a$ and $\gamma(1) = b$.
\end{lemma}
%
%

\begin{lemma}
\label{lem:arctopoly}
Let $\gamma$ be a simple path in an open
set $U \subseteq \mathbb{R}^2$,
where $\gamma(0) = a$ and $\gamma(1) = b$.
Choose an integer $n > 0$.
Suppose that the diameter of $\gamma$ is 
less than $\epsilon$.
Then there exists a simple polygonal path in $U$ 
of diameter less than 
$\left(1 + \frac{1}{2^{n-1}}\right) \epsilon$ 
which starts at
$a$ and ends at $b$.
\end{lemma}
%

\begin{remark}
\label{rem:ulcpoly}
\textnormal{
Note that we are not attempting to approximate
$\gamma$ arbitrarily closely by a simple polygonal path.
We are merely trying to find a polygonal path in $U$
connecting $a$ and $b$ given an arbitrary path in $U$
connecting these two points.
This way,
we may use polygonal paths instead of arbitrary paths.
In particular,
if a uniformly locally arcwise connected open set 
contains a path
of diameter less than $\epsilon$ connecting $z_1$
and $z_2$ if $|z_1 - z_2| < \epsilon$,
then it also contains a simple polygonal path of
diameter less than $2 \epsilon$ connecting these
points.
Hence,
without loss of generality,
we may use simple polygonal paths in place of paths
in arguments involving a uniformly locally 
arcwise connected
open set.
}
\end{remark}

\begin{lemma}
\label{lem:tubepath}
Let $U$ be an open
subset of $\mathbb{R}^2$ and let $\gamma \subset U$
be a simple polygonal path with $\gamma(0) = a$ and
$\gamma(1) = b$ of diameter less than $\epsilon$.
Choose $\delta_0 > 0$.
Let $c = \gamma(t_c)$ for some $t_c \in (0, 1)$ and 
suppose that $\zeta \in B(c, \delta) \backslash \gamma$
where $0 < \delta < \delta_0$ is chosen so that
$B(c, \delta) \subset U$.
Then there exists a simple polygonal path 
$\tilde{\gamma}
\subset U$ of diameter less than 
$\epsilon + 2 \delta_0$
such that $\tilde{\gamma}(0) = \zeta$,
$\tilde{\gamma}(1) = b$, and $\tilde{\gamma} \cap
\gamma = \{b\}$.
\end{lemma}
\textit{Proof.}
While it is possible to give an explicit construction
of the desired path,
we will spare the reader the details and instead
follow a suggestion of D. Sarason.

Cover each $z \in \gamma$ with a ball of radius
$r_z \leq \delta$ such that $B(z, r_z) \subset U$.
Include $B(c, \delta)$ to cover $\{\zeta\}$.
Because $\gamma \cup \{\zeta\}$ is compact,
we may choose a finite subcover,
which consists of balls centered at points
$z_1, z_2, ..., z_N$ of $\gamma$.
Let $G = \bigcup_{j=1}^N B(z_j, r_{z_j})$.
By construction,
the diameter of $G$ is at most 
$\epsilon +2 \delta$.
Let $B = G \backslash \gamma$.
Clearly,
$B$ is open.
It also follows that $B$ is connected.
(Let $F_1 = \mathbb{R}^2 \backslash G$ and
$F_2 = \gamma$.
Then $F_1 \cap F_2 = \varnothing$.
Since $G$ is open,
$F_1$ is closed.
By construction, 
$\mathbb{R}^2 \backslash F_1 = G$ is connected.
Clearly,
$F_2 = \gamma$ is closed.
Since $\gamma$ is a simple arc,
it does not separate the plane and 
$\mathbb{R}^2 \backslash F_2$ is connected.
It follows from Corollary 1, page 112 of 
\cite{new:plane}
that $B = \mathbb{R}^2 \backslash (F_1 \cup F_2)$
is connected.)

Choose $r$ such that $B(b, r) \subset B(z_N, r_N)$
and such that $B(b, r)$ intersects exactly one
line segment of $\gamma$
(this is possible because $b$ is the last point of
$\gamma$). 
Choose $\xi \in B(b, r)$.
It is well known (see,
for example,
page 15 of \cite{con:cx1})
that an open set $U \subset \mathbb{R}^2$ is
connected if and only if every pair of
points in $U$ can be connected by a polygonal path
contained in $U$.
In particular,
we can find a polygonal path,
say $\Gamma$, 
in $B$ connecting
$\zeta$ and $\xi$.
By construction,
$\gamma \cap \Gamma = \varnothing$.
By our choice of $r$,
we may find a line segment $\ell$ in 
$B(b, r) \subset B(z_N, r_N)$
connecting $\xi$ to $b$, 
such that $\gamma \cap \ell = \{b\}$.
Apply Lemma \ref{lem:simplepoly} 
to find a simple polygonal path
$\tilde{\gamma} \subseteq \Gamma \cup \ell$
connecting $\zeta$ and $b$.
By construction,
the diameter of 
$\tilde{\gamma} \leq \epsilon + 2 \delta
< \epsilon + 2 \delta_0$.
By construction,
$\tilde{\gamma}$ is a simple polygonal path,
starting at $\zeta$ and ending at $b$,
such that $\tilde{\gamma} \cap \gamma = \{b\}$.
\qed

We now apply these tools to build an asymptotic path.

\begin{lemma}
\label{lem:endcut}
Suppose that $U$ is a connected, open subset of 
$\mathbb{R}^2$ and
that $\{\zeta_n\} \subset U$ such that
$\zeta_n \rightarrow \zeta_0 \in \partial U$,
where $\zeta_0$ is finite.
Also suppose that there exists $\epsilon > 0$ such that
$B(\zeta_0, \epsilon) \cap U$ is uniformly 
locally arcwise connected.
Then there exists an asymptotic path
$\Gamma \subset U$ containing a subsequence of
$\{\zeta_n\}$ such that $\Gamma \rightarrow \zeta_0$. 
\end{lemma}
\noindent
\textit{Proof.}
Without loss of generality,
we will let $U = U \cap B(\zeta_0, \epsilon)$.
Choose $\epsilon_0 > 0$ such that 
$\epsilon_0 < \epsilon$.
Since U is ulac
and by Remark \ref{rem:ulcpoly},
there exists $\delta_0 > 0$ such that whenever
$z_1, z_2 \in U$ with $|z_1 - z_2| < \delta_0$,
then there exists a simple polygonal path in 
$U$ connecting
$z_1$ and $z_2$ of diameter less than $\epsilon_0 / 2$.

By passing to a subsequence,
we may assume that 
$|\zeta_1 - \zeta_0| < \delta_0 / 2$
and that $|\zeta_{n+1} - \zeta_0| <
\frac{1}{2} |\zeta_n - \zeta_0|$.
Thus, 
$|\zeta_j - \zeta_0| < \delta_0 / 2^j$
and $|\zeta_j - \zeta_k| 
\leq |\zeta_j - \zeta_0| + |\zeta_k - \zeta_0|
< \delta_0 / 2^{j-1}$ for 
$k > j$.
We will construct an asymptotic path
which consists of a countable number of
simple polygonal paths in $U$
and which
contains $\{w_m\}_{m=1}^\infty$, 
a subsequence of $\{\zeta_n\}$.

Let $w_0 = \zeta_0$ and $w_1 = \zeta_1$.
Let $\rho_1 = \epsilon_0$ and let $d_2 = \rho_1 / 4$.
Choose $j > 1$ so that 
(i) there exists a path in $U$
of diameter less that $d_2$
connecting $z$ and $\xi$ whenever
$|z - \xi| < \delta_0 / 2^{j-1}$
(by ulac of $U$) and 
(ii) $d(\zeta_j, w_0) < \frac{1}{4}\, d(w_1, w_0)$
(since $d(\zeta_n, w_0)$ is strictly
decreasing).
Let $w_2 = \zeta_j$.
Since $|w_2 - w_1| = |\zeta_j - \zeta_1| < \delta_0$,
there is a simple polygonal path,
say $\gamma_{1,2}$,
in $U$ of diameter less than 
$\epsilon_0 / 2$ connecting $w_1$ and $w_2$.
Let $\Gamma_1= \varnothing$.
Let $\Gamma_2 = \gamma_{1,2}$\ and
let $\rho_2 = d(\gamma_{1,2}, w_0) = d(\Gamma_2, w_0)$.

We will build the other segments inductively
for $n > 1$.
Let $\gamma_{n, n+1}$ be the simple polygonal path
constructed in this procedure connecting
$w_{n}$ and $w_{n+1}$.
Let 
$\Gamma_{n+1} = \Gamma_n\, \cup\, \gamma_{n, n+1}$.
Let $\rho_{n+1} = d(\Gamma_{n+1},\, w_0)$.
The path $\gamma_{n, n+1}$ is chosen so that
$\gamma_{n, n+1} \cap \Gamma_{n} = \{w_n\}$
and 
$\gamma_{n, n+1} \cap \Gamma_{n-1} = \varnothing$.
We also require that $\gamma_{n, n+1} \subset
U \cap B(w_0, \frac{7}{8}\, \rho_{n-1})$;
this guarantees that $\rho_n \rightarrow 0$.
Let
$\Gamma = \lim_{n \rightarrow \infty} \Gamma_n$.
If each $\gamma_{n, n+1}$ for $n > 1$ satisfies these
assumptions,
then $\Gamma$
is a simple asymptotic path which approaches $\zeta_0$
and contains a subsequence of $\{\zeta_j\}$.

It remains to construct $\gamma_{n, n+1}$ for $n > 1$.
It is enough to choose $w_{n+1}$ so that 
(i) $\gamma_{n, n+1}\,
\cap\, \Gamma_{n} = \{w_n\}$, 
(ii) 
$\gamma_{n, n+1} \cap\, \Gamma_{n-1} = \varnothing$
and (iii)
$\gamma_{n, n+1} \subset U \cap B(w_0, \frac{7}{8}\, 
\rho_{n-1})$.

We start by choosing $j$ so that
(a) $|\zeta_j - w_0| < \rho_n / 4$ and
(b) for all $z \in U$ such that $|z - \zeta_j| <
\delta_0 / 2^{j - 1}$,
there exists a simple polygonal path in $U$ from
$z$ to $\zeta_j$ of diameter less than 
$d_{n+1} = \min(\epsilon_0 / 2^{n+1}, \rho_n / 4)$.
Condition (a) forces $\zeta_j$ to be closer to $w_0$
than $\Gamma_n$ is.
Condition (b) lets us control the diameter of
simple polygonal paths from
from $\zeta_j$ to $\zeta_k$ for all $k > j$.
It is clear that (a) is possible and that
$j$ will be greater than the corresponding $k$
in $w_j=\zeta_k$
for $w_1, w_2, ..., w_n$.
Remark \ref{rem:ulcpoly} allows us to satisfy
(b).
Let $w_{n+1} = \zeta_j$.

Let $\eta$ be a simple polygonal path in $U$ connecting
$w_n$ and $w_{n+1}$ of diameter less than $d_n$
(this is possible by the previous stage of the 
construction).
Recall that $d_n \leq \rho_{n-1} / 4$ and that
$\rho_n \leq d(w_n, w_0) < \rho_{n-1} / 4$.
Recall also that $\Gamma_1 = \varnothing$,
so we only need to check (ii) for $n > 2$.
When $n > 2$,
$\eta \subset B(w_0, r)$ where
$r \leq |w_n - w_0| + d_n < \rho_{n-1}/4 + \rho_{n-1}/4
= \rho_{n-1} / 2$,
so $\eta$ will not intersect
$\Gamma_{n-1}$.

The only possible way for $\eta \cup \Gamma_n$
to fail being a simple polygonal path is if
$\eta$ intersects $\gamma_{n-1, n}$ at a point
other than $w_n$.
Let $c$ be the first point of intersection
of $\eta$ with $\gamma_{n-1, n}$ 
as we travel from $w_{n+1}$ to $w_n$.
Choose $\delta_0$ (in Lemma \ref{lem:tubepath})
less than $ d_{n+1}$.
Choose $\zeta$ in Lemma \ref{lem:tubepath} so that
$\zeta$ is in the portion of $\eta$
starting from $w_{n+1}$ before $c$.
Apply Lemma \ref{lem:tubepath} to
construct $\tilde{\eta}$ connecting $\zeta$ to $w_n$
and
let $\eta^\prime$ be the portion of $\eta$ from
$w_{n+1}$ to $\zeta$.
Let $\tilde{\gamma}_{n,n+1} = \eta^\prime \cup 
\tilde{\eta}$ and apply Lemma \ref{lem:simplepoly} to 
find $\gamma_{n, n+1}$.
This new path satisfies (i) above.
We note that if (iii) holds, so will (ii).
We now check (iii):
By Lemma \ref{lem:tubepath},
$\gamma_{n, n+1} \subset U \cap B(w_0, r_n)$
where
$r_n < \rho_{n-1}/2 + (d_n + 2 \delta_0)$
$< \rho_{n-1} / 2 + \rho_{n-1} / 4 + 2 \rho_n / 4$
$< 3 \rho_{n-1} / 4 + \rho_{n-1} / 8$
$= \frac{7}{8}\, \rho_{n-1}$.
This completes the construction.
\qed

\begin{remark}
\textnormal
{
It is known that a boundary point
of a uniformly locally connected open set is 
accessible.
Note that a ulac open set is also uniformly 
locally connected,
so a boundary point of a ulac set will also be
accessible.
The difficulty here is that we want the end-cut
to contain a subsequence of $\{f(z_n)\}$ where
$\{z_n\} \rightarrow \partial R \cup \{\infty\}$
and $f(z_n) \rightarrow w_0$.
}
\end{remark}

\begin{lemma}
\label{lem:imagepath}
Suppose that $U$ is an open subset of $\mathbb{R}^2$ and
that $\{\zeta_n\} \subset U$ such that
$\zeta_n \rightarrow \zeta_0 \in \partial U$,
where $\zeta_0$ is finite.
Also suppose that there exists $\epsilon > 0$ such that
$B(\zeta_0, \epsilon) \cap U$ consists of 
finitely many components and that
each component is uniformly locally arcwise connected.
Then there exists a simple continuous curve 
$\Gamma \subset U$ containing a subsequence of
$\{\zeta_n\}$ such that $\Gamma \rightarrow \zeta_0$. 
\end{lemma}
\noindent
\textit{Proof.}
Without loss of generality, 
we assume that our sequence consists of distinct 
points, 
that $\{\zeta_n\} \subset B(w_0, \epsilon)$,
and that $|\zeta_{n+1} - \zeta_n|< \frac{1}{2}
|\zeta_n - \zeta_{n-1}|$.
By the pigeonhole principle,
we may find a subsequence of $\{\zeta_n\}$
that lies in one of the components of
$B(\zeta_0, \epsilon) \cap U$,
say $\Omega_0$.
Our result follows by applying Lemma \ref{lem:endcut}
to $\Omega_0$.
\qed

We now apply Theorem \ref{thm:part} to lift our
asymptotic path.

\begin{lemma}
\label{lem:dompath}
Let $R \subseteq \mathbb{R}^2$ be open and
let $f:R \rightarrow \mathbb{R}^2$ be $C^1$.
Let $w_0 \in C(f)$ and $\{z_n\} \subset R$
be such that $f(z_n) \rightarrow w_0$ and
$z_n \rightarrow z_0 \in \partial R \cup \{\infty\}$.
Suppose that $\{z_n\} \subset R_0$ where 
$R_0$ is a component
of $R \backslash f^{-1}(f(S) \cup C(f))$.
Let $w_n = f(z_n)$ and let  $\Omega_0 = f(R_0)$.
Suppose also that there exists $\epsilon > 0$
such that $B(w_0, \epsilon) \cap \Omega_0$ can
be partitioned into a finite number of components,
each of which is uniformly locally arcwise connected.
Then there exists a simple, half-open path 
$\gamma \subset R_0$ which
contains a subsequence of the $\{z_n\}$.
\end{lemma}
\noindent
\textit{Proof.}
By Theorem \ref{thm:part},
$\Omega_0$ is a component of 
$\mathbb{R}^2 \backslash (f(S) \cup C(f, \infty))$
and $f|_{R_0}$ is a covering map.
Let $N = Val(f|_{R_0})$.
By the pigeonhole principle, a subsequence of
$\{w_n\}$
is contained in one of the components of
$B(w_0) \cap \Omega_0$.
By Lemma \ref{lem:imagepath},
there exists a simple asymptotic path
$\Gamma \subset \Omega_0$
that contains a pairwise distinct
subsequence of $\{w_n\}$ such
that 
$\Gamma \rightarrow w_0$.
We will denote this subsequence by $\{w_n\}$ and
the corresponding subsequence of $\{z_n\}$ by $\{z_n\}$.
Because $\Gamma$ is the continuous image of 
$[0, \infty)$
(whose fundamental group is trivial) and 
$f|_{R_0}$ is a covering map,
we can lift $\Gamma$ from $\Omega_0$
to $R_0$.
Since $f$ is locally 1-1 in $R_0$
and $\Gamma$ is simple,
$\Gamma$ lifts to $N$ distinct,
simple half-open paths in $R_0$.
By the pigeonhole principle,
at least one of these paths,
say $\gamma$, 
contains an infinite subsequence of the $\{z_n\}$.
\qed

\begin{lemma}
\label{lem:asympath}
Let $R = \mathbb{R}^2$.
Let 
$f$, $w_0 \in C(f, \infty)$ and $\{z_n\}$ 
be defined as in Lemma \ref{lem:dompath}.
Suppose that $\{z_n\} \subset R_0$ where $R_0$ is a 
component of $\mathbb{R}^2 \backslash 
f^{-1}(f(S) \cup C(f, \infty))$.
Let $w_n = f(z_n)$ and let  $\Omega_0 = f(R_0)$.
Suppose also that there exists $\epsilon > 0$
such that $B(w_0, \epsilon) \cap \Omega_0$ can
be partitioned into a finite number of components,
each of which is uniformly locally arcwise connected.
If either (i) $Val(f, w_0) < \infty$ or 
(ii) $Val(f, w_0) = \infty$ and $w_0 \notin f(S)$,
then $w_0$ is an asymptotic value of $f$.
\end{lemma}
\noindent
\textit{Proof.}
By Lemma \ref{lem:dompath},
there exists a simple, half-open path
$\gamma \subset R_0$ such that
$\gamma$ contains a subsequence of $\{z_n\}$ and
$f(\gamma) \rightarrow w_0$.
It remains to show that $\gamma \rightarrow \infty$.
Suppose not.
Then there exists $\rho > 0$ such that every unbounded 
subarc of $\gamma$ intersects $B(0, \rho)$.
In other words,
we can find a subsequence of $\{z_n\}$
such that $|z_n|$ is strictly increasing and
such that $\gamma_n$, 
the lift in $\gamma$ starting at $z_n$ of the 
path joining $w_n$ and $w_{n+1}$ in $f(\gamma)$,
intersects $B(0, \rho)$.
In particular, 
choose an integer $m > \rho$.
Then there exists $N$ such that
$|z_n| > m$ for all $n > N$ in our subsequence.
Let $C_m = \{z : |z| = m\}$.
Then,
for each $n > N$ in our subsequence,
there exists $\zeta_n \in \gamma_n \cap C_m 
\subset \gamma \cap C_m$. 
By construction,
the $\zeta_n$ will be pairwise distinct.
Since $f(\gamma) \rightarrow w_0$,
$f(\zeta_n) \rightarrow w_0$.
Since $C_m$ is compact,
$\{\zeta_n\}$ has a convergent subsequence.
Hence, $\zeta^{(m)}$, the limit of this subsequence,
is a preimage of $w_0$ in $C_m$.
We can repeat this construction for each integer
$m > \rho$.
Since the $C_m$ are pairwise disjoint,
the $\zeta^{(m)}$ are distinct and
$Val(f, w_0) = \infty$, 
a contradiction when $Val(f, w_0)$ is finite.

On the other hand,
if $Val(f, w_0) = \infty$,
then we may assume that $w_0 \notin f(S)$.
This will also lead to a contradiction.
Consider the closed annulus 
$A = \{z: m \leq |z| \leq m + 1\}$,
where $m > \rho$.
By the construction above,
each circle $C_r$ where $m \leq r \leq m+1$
contains a point $\zeta_r$ such that $f(\zeta_r) = w_0$.
Thus the $\zeta_r$ accumulate in $A$,
a compact set.
Let $\zeta \in A$ be one such accumulation point.
Then $f$ cannot be $1-1$ in any neighborhood of
$\zeta$.
By the inverse function theorem, 
$J_f(\zeta) = 0$.
Hence $\zeta\in S$ and $w_0 = f(\zeta) \in f(S)$,
a contradiction.
\qed

\begin{remark}
\textnormal{
The construction used in the proof of 
Lemma \ref{lem:asympath} uses our additional 
requirement that our lift is a simple path,
which guarantees that each time the path
intersects the circle $C_m$, 
it does so at different points.
Thus $\{\zeta_n\}$ consists of an infinite number
of distinct points.
}
\end{remark}

\begin{remark}
\textnormal{
The construction used in the proof of 
Lemma \ref{lem:asympath} can be extended
to the case where $R$ is an open disk,
$w_0 \in C(f)$, and 
$R_0$ is a component of 
$R \backslash (f(S) \cup C(f))$ by
intersecting $R$ with a sequence of
circles that approach $\partial R$.
The lift of $\gamma$ approaches $\partial R$;
otherwise,
we can argue as above
to get the contradiction that $Val(f, w_0) = \infty$.
However,
it is unclear whether the lift approaches a
particular point on $\partial R$,
so we cannot claim that $w_0$ is an asymptotic value.
}
\end{remark}

\begin{lemma}
\label{lem:altfindseq}
Let $R \subseteq \mathbb{R}^2$ be open and
let $f:R \rightarrow \mathbb{R}^2$ be $C^1$.
Suppose that $w_0 \in C(f)$ and that $S$ is 
nowhere dense in $R$.
Suppose that there exists $\epsilon > 0$
such that 
$B(w_0, \epsilon) \backslash (f(S) \cup C(f))$ 
consists of finitely many components 
and such that $C(f)$ is nowhere dense in
$B(w_0, \epsilon)$.
Then there exists a sequence $\{\zeta_n\} \subset R_0$,
where $R_0$ is a component of
$R \backslash f^{-1}(f(S) \cup C(f))$,
and
such that $\zeta_n \rightarrow \zeta_0 \in \partial R\,
\cup\, \{\infty\}$.
\end{lemma}
\noindent
\textit{Proof.}
By definition,
there exists $\{z_n\} \subset R$ such that
$w_n = f(z_n) \rightarrow w_0$ and 
$z_n \rightarrow z_0 \in
\partial R
\cup \{\infty\}$.
To construct $\{\zeta_n\}$,
we start by guaranteeing that 
this sequence does not intersect the partitioning
set of the preimage.
Without loss of generality, 
we may choose a subsequence of $\{w_n\}$ such that 
$|w_{n+1} - w_0| < \frac{1}{2} |w_n - w_0|$
for each $n$.
We may choose $N$ such that for $n > N$,
$w_n \in B(w_0, \epsilon)$.
Let $\epsilon_n = \frac{1}{2}
\min(\epsilon - |w_n - w_0|, |w_n - w_0|)$.
Choose $0 < \delta < 1/n$ such that
$f(B(z_n, \delta)) \subseteq B(w_n, \epsilon_n)$.
Since $S$ is nowhere dense, 
$\exists z'_n \in B(z_n, \delta) \backslash S$ and
we may choose an open set $U \subseteq B(z_n, \delta)$ such that
$z'_n \in U$ and $f$ is a homeomorphism on $U$.
Recall that $f(S) \cup C(f)$ is closed.
By Sard's theorem
(see,
for example,
page 69 of \cite{hir:difftop}),
$f(S)$ has empty interior.
Thus,
the restriction of
$f(S) \cup C(f)$ to $B(w_0, \epsilon_0)$
has empty interior iff the restriction of $C(f)$ 
does
(see Remark 3.8 of \cite{neu:valence}).
In particular,
since the closed set
$C(f)$ is nowhere dense in $B(w_0, \epsilon_0)$,
$f(S) \cup C(f)$ is nowhere dense in 
$B(w_0, \epsilon_0)$.
Therefore there exists 
$w'_n \in f(U) \subseteq B(w_n, \epsilon_n)$
such that $w'_n \notin f(S) \cup C(f)$.
Let $\zeta_n = f^{-1}(w'_n) \cap U$.
Then
$\zeta_n \rightarrow z_0$ and 
$f(\zeta_n) \rightarrow w_0$.

We then show that there is a subsequence of 
$\{\zeta_n\}$ in one component of the preimage.
Note that a component of 
$B(w_0, \epsilon) \backslash (f(S) \cup C(f))$
is a subset of some component of
$\mathbb{R}^2 \backslash (f(S) \cup C(f))$.
Since
$B(w_0, \epsilon) \backslash (f(S) \cup C(f))$
consists of finitely many components,
$B(w_0, \epsilon)$ intersects finitely many
components of the partition of the image.
Since each region in the partition of the image has
finite valence,
$\{\zeta_n\}$ is contained in finitely many regions
of the partition of the preimage.
Applying the pigeonhole principle gives us a 
subsequence of
$\{\zeta_n\}$ contained in one component of the 
partition of the preimage.
\qed

\noindent
\textbf{Proof of Theorem \ref{thm:locconpath}.}
By Lemma \ref{lem:altfindseq},
there exists a sequence $\{\zeta_n\} \subset R_0$,
where $R_0$ is a component of 
$\mathbb{R}^2 \backslash f^{-1}(f(S) \cup
C(f, \infty))$, such that
$\zeta_n \rightarrow \infty$ 
and $f(\zeta_n) \rightarrow w_0$.
The result then follows from Lemma \ref{lem:asympath}.
\qed

\section{Examples and discussion}

\begin{example}
\label{ex:nono}
A harmonic function with
a cluster set with non-empty interior such that
$f(\mathbb{C}) \backslash (f(S) \cup C(f, \infty))
\neq \varnothing$ and such that
$\overline{f(S)} \backslash f(S) \neq \varnothing$.
\begin{displaymath}
f(z) = \Re\, e^z + \frac{i}{2}\, \Im\, z^2
= e^x \,\cos\, y + i\, x y
\end{displaymath}
\end{example}
A calculation shows that 
\begin{displaymath}
S = \{z = x + i\, y: x = -y\, \tan\, y,\ 
y \neq (2 k + 1) \pi /2,  k \in \mathbb{Z}\}
\end{displaymath}
\begin{displaymath}
f(S) = \{e^{-y\, \tan\,y}\, \cos\, y -
i\, y^2\, \tan\,y,\ 
y \neq (2 k + 1) \pi /2,  k \in \mathbb{Z}\}
\end{displaymath}
$f(S)$ intersects the real axis only at 
$\Re\,w = \pm 1$.
These are the only points of intersection of
$f(S)$ with the vertical lines $\Re{w} = \pm 1$.
This follows by first noting that
$e^{-y\, \tan\, y}\, \cos\,y$ is
strictly monotone on the interval 
$(2 k - 1) \pi / 2 < y < (2 k + 1) \pi /2$
when $k \neq 0$.
When $k = 0$,
$e^{y\, \tan\, y} \geq 1$ and thus
will intersect $\cos\,y$ at most once on our interval.
By considering $y = k \pi - b / (k \pi)^2$,
we can check that $f(S)$ is arbitrarily close to
$(-1)^k + i\,b$ as $|k| \rightarrow \infty$.
Thus,
$\{w: \Re\,w = \pm 1, \Im\, w \neq 0\} 
\subseteq \overline{f(S)} \backslash f(S)$.

We claim that $C(f, \infty) = 
\{w: -1 \leq \Re\,w \leq 1\}$.
Let $w = a + ib$. 
Suppose that $f(z_n) \rightarrow w$ as 
$z \rightarrow \infty$. 
If $a > 1$,
we must have $\lim x_n > 0$.
If $x_n \rightarrow \infty$,
we must have $y_n \rightarrow (2 k + 1) \pi /2$;
however, 
$x_n y_n \rightarrow b$ is then impossible.
Moreover,
if $y_n \rightarrow \infty$,
we must have $x_n \rightarrow 0$. 
Thus, 
$\Re w = a > 1$ is omitted from $C(f, \infty)$.
Similarly,
$\Re\,w < -1$ is omitted from $C(f, \infty)$.
Thus,
if $w \in C(f, \infty)$,
then $|\Re\,w| \leq 1$.
Now consider the case when $|a| \leq 1$.
Choose $b \in \mathbb{R}$.
It's clear that $\cos\,y$ and
$a\,e^{-b/y}$ have an infinite number of intersections
when 
(i) $|y| > 1$ and, 
(ii) when $b$ is nonzero,
$y$ also must have the same sign as $b$.
Let $\{y_n\}$ be chosen from these points of 
intersection and let 
$\{z_n = (b / y_n) + i y_n\}$.
Thus,
$Val(f, w) = \infty$ for $w \in C(f,\infty)$.

We also note that 
$f(\mathbb{C}) \backslash (f(S) \cup C(f, \infty)) \neq 
\varnothing$.
For example,
$w = a > 1$, 
so $w \notin f(S) \cup C(f, \infty)$.
Then $f^{-1}(w) = \{\log\,a\}$.
Note that $f$ omits the real axis for $a < -1$.

\begin{example}
An example for Theorem \ref{thm:clusvals}.
\begin{displaymath}
f(z) = e^z + i\, \Im z = 
e^x \cos y + i\, (e^x \sin\, y + y)
\end{displaymath}
\end{example}
A calculation shows that $\{z_k = 
\log ((4 k + 3) \pi / 2) + i\, (4 k + 3) \pi / 2:
k \in \mathbb{Z}^+\} \subset f^{-1}(0)$.
This sequence lies in the right half plane,
with $\Re{z_k} > 1$ and $z_k \rightarrow \infty$.
We check the conditions of Theorem \ref{thm:clusvals}.
We note that no point in the critical set of $f$ is
in the right half plane,
because $z = x + iy \in S$ implies that 
$\cos y = -e^x$.
Since $\{z_k\}$ is in the right half plane,
it is not
in $S$ and condition (1) is satisfied.
Further,
$J_f(z_k) = (\frac{4k + 3) \pi}{2})^2$,
so (2) is satisfied.
Let
$B_k = \overline{B(z_k, \log(2)/2)}$.
For $z \in B_k$,
$\Re{z} \leq  \Re z_k + \log(2) / 2 = 
\log(\frac{(4k+3) \pi}{\sqrt{2}})$
and
$\Re{z} \geq \Re z_k - \log(2)/2 = 
\log(\frac{(4k+3) \pi}{2 \sqrt{2}})
> \log(\frac{3 \pi}{2 \sqrt{2}})$.
A calculation shows that condition (3) is satisfied,
since $J_f(z) / J_f(z_k) <
2 + \frac{2 \sqrt{2}}{3 \pi}$ on each $B_k$.
We can thus apply Theorem \ref{thm:clusvals}
if
$|\mu_f(z)| = |f_{\overline{z}}(z) / f_z(z)|$
is bounded away from 1 on each $B_k$.
For $z \in B_k$,
$|\mu_f(z)| = 1 / |2 e^z + 1| \leq
1 / (|2 e^z| - 1) = 1 / (2 e^{\Re(z)} - 1)
< 1 / (\frac{3 \pi}{\sqrt{2}} - 1) < \frac{1}{2}$.
So $\{z_k\}$ satisfies the conditions
of Theorem \ref{thm:clusvals}.
Thus,
$f$ has infinite valence in a neighborhood
of $w_0 = 0$ and the interior of
$C(f, \infty)$ in non-empty.
From a result of M. Balk
(\cite{ba:polyentire}, page 107),
we know that $C(f, \infty) = \mathbb{C}$ and
that every $w \in \mathbb{C}$ has infinite valence.

\begin{example}
The cluster points of $f(z) = z + \Re\,e^z$
are all asymptotic values.
\end{example}
$C(f, \infty)$ consists of the horizontal lines
where the imaginary part of $w$ is an odd multiple
of $\pi/2$.
Note that while the horizontal lines in the 
cluster set can be thought
of as horizontal asymptotes to the image of the
critical set as $z \rightarrow \infty$,
$f(S) \cap C(f, \infty) = \varnothing$.
Further,
$Val(f)$ is finite.
For details,
see Example 4.2 in \cite{neu:valence}.
If we choose $w \in C(f, \infty)$,
the conditions in Theorem \ref{thm:locconpath}
are satisfied and $w$ is an asymptotic value.

\begin{example}
The cluster points of $f(z) = 2 (\Re(z^3) + iz)$
are all asymptotic values.
\end{example}
$C(f, \infty)$ is the real axis
and $Val(f)$ is finite.
For details,
see Example 4.3 in \cite{neu:valence}.
While $f(S) \cap C(f, \infty) = \varnothing$,
$f(S)$ is asymptotic to $C(f, \infty)$ as
$z \rightarrow \infty$.
If we choose $w \in C(f, \infty)$,
the conditions in Theorem \ref{thm:locconpath}
are satisfied and $w$ is an asymptotic value.

\begin{example}
The finite valence cluster points of 
$f(z) = 2 [\Re(z^2) + i(\Re(z) - \Im(z))]$
are all asymptotic values.
\end{example}
$C(f, \infty)$ is the real axis.
Note that every cluster point other than the
origin is omitted and
that the origin has infinite valence.
Also note that $f(S) = \{0\}$.
For details,
see Example 4.4 in \cite{neu:valence}.
If we choose $w \in C(f, \infty) \backslash \{0\}$,
the conditions in Theorem \ref{thm:locconpath}
are satisfied and $w$ is an asymptotic value.
Note that $w = 0$ is also an
asymptotic value;
choose $\Re z = \Im z$ and
let $z \rightarrow \infty$.

\begin{example}
The finite valence cluster points of 
$f(z) = \Im(\exp(z^2)) + i \Im({z^2})$
are all asymptotic values.
\end{example}
Let $z = x + iy$.
Calculations show that
\[
S = \{xy = 0\} 
\cup \{y = k \pi / (2x): k \in \mathbb{Z}\}
\]
\[
f(S) = \{i k \pi: k \in \mathbb{Z}\}
\]
One can easily check that $Val(f, w) = \infty$
for $w \in f(S)$
and that each such $w$ is an asymptotic value.
Thus $f(S) \subseteq C(f, \infty)$.
Moreover,
points other than $w = i k \pi$
on the imaginary axis and
on the horizontal lines $\Im w = k \pi$
are omitted.
The valence is constant within each of the
remaining horizontal strips in the right half plane
and the left half plane.
When the points in a horizontal strip in the right
half plane have valence 2,
points in that strip in the left half plane are
omitted.
Similarly,
if the points in a horizontal strip in the right
half plane are omitted,
points in that strip in the left half plane have
valence 2.
Whether the valence is 2 or 0 depends on
whether or not
$\Re w$ and $\sin(\Im w)$ have the same sign.
We note that $S$ is nowhere dense.
We see from Theorem \ref{thm:part} that
$C(f, \infty)$ contains the horizontal lines
$\Im w = k \pi$,
as well as the imaginary axis,
because none of these points has a 
neighborhood where the valence is constant.
In fact, 
this exhausts $C(f, \infty)$.
(Otherwise,
we would have $\{z_n\} \rightarrow \infty$
with $f(z_n) \rightarrow w \in C(f, \infty)$,
satisfying 
$\Re f(z_n) \rightarrow \Re w \neq 0$
and
$2 x_n y_n \rightarrow \Im w \neq k \pi$.
Since this implies that
$\sin(2 x_n y_n) \rightarrow \sin(\Im w) \neq 0$,
we see that $\Re f(z_n)$ having a finite limit requires 
$|y_n| \rightarrow \infty$
and thus $\Re w$ must be zero,
a contradiction.)
If $w \in C(f, \infty) \backslash f(S)$,
then
the conditions in Theorem \ref{thm:locconpath}
are satisfied and $w$ is an asymptotic value.
For example,
if $w = ib$ for $b \neq k \pi$,
let $y = b / (2x)$ and let $x \rightarrow 0^+$.
As another example,
if $w = 1$,
let $y = \exp(-x^2) / (2x)$ and let
$x \rightarrow \infty$.

\begin{remark}
\textnormal{
Example \ref{ex:nono} above illustrates why we require
that the cluster set have empty interior
in some neighborhood of the cluster point.
Choose $w \in C(f, \infty)$ such that $\Re w = \pm 1$.
Although for $\epsilon > 0$,
$B(w, \epsilon) \backslash C(f, \infty) 
\neq \varnothing$,
$w$ is not an asymptotic value of $f$.
In order for $\Re f$ to approach $\Re w$,
our asymptotic path would 
either have $\Re z$ bounded
(but then $\Im z \rightarrow \pm \infty$,
which gives a contradiction 
because of the oscillation in $\Re f$)
or have $\Re z \rightarrow +\infty$
(however,
$\Im f \rightarrow \Im w$ requires that 
$\Im z \rightarrow 0$,
which is impossible as $|\Re f|$ would then be
unbounded).
}
\end{remark}

\begin{corollary}
\label{cor:arbclose}
Suppose that $f: \mathbb{C} \rightarrow \mathbb{C}$
is harmonic and
that $Val(f, w)$ is finite
for all $w \in \mathbb{C}$.
Suppose that $S$ is nowhere dense.
Let $w_0 \in C(f, \infty)$.
Then at least one of the following holds
for some subsequence of
any given sequence $\{z_n\}$
such that $z_n \rightarrow \infty$ and
$f(z_n) \rightarrow w_0$:
\begin{enumerate}
\item
$d(z_n, S) \rightarrow 0$ as
$z_n \rightarrow \infty$.
\item
$J_f(z_n) \rightarrow 0$ as
$z_n \rightarrow \infty$.
\item
There is no $\rho > 0$ such that
$|J_f(z) / J_f(z_n)| < M$ on $B(z_n, \rho)$ for all $n$
for some finite constant
$M$.
\item
There is no $\rho > 0$ such that
$sup_n\{|\mu_f(z)|: z \in \overline{B(z_n, \rho)}\}
< 1$.
\end{enumerate}
\end{corollary}
\noindent
\textit{Proof.}
By Theorem 3.9 in \cite{neu:valence},
if $S$ is nowhere dense and all points have finite 
valence,
then the cluster set has empty interior.
Then apply Theorem \ref{thm:clusvals}.
\qed

\begin{remark}
\textnormal{
If $f$ is an entire harmonic function,
with finite valence and $S$ nowhere dense,
Corollary \ref{cor:arbclose} holds.
Let $w_0$ be a cluster point of $f$.
If condition (1) 
holds,
$\{z_n\}$ gets arbitrarily close to $S$
(though this does not guarantee that 
$J_f(z_n) \rightarrow 0$).
If condition (2) holds,
there is no region containing the
$\{z_n\}$ where the function is
quasiconformal or quasiregular.
}
\end{remark}

\begin{remark}
\textnormal{
Comparing the Conjecture above with 
Theorem \ref{thm:locconpath},
it remains to show that harmonicity implies
condition (2) on $B(w_0, \epsilon)$ in 
Theorem \ref{thm:locconpath}.
It would also be interesting to find a 
necessary and sufficient condition for
$C(f, \infty)$ to have empty interior.
}
\end{remark}

\section*{Acknowledgements}
The author is grateful to D. Sarason for reading
earlier versions of this write-up
and for his helpful comments and suggestions.
The author thanks P. Poggi-Corradini and
D. Auckly for interesting discussions.

\bibliographystyle{amsplain}

\end{document}